\documentclass[oneside,english]{article}
\usepackage[T1]{fontenc}
\usepackage[latin9]{inputenc}
\usepackage{amsmath}
\usepackage{amsthm}
\usepackage{amsfonts}
\usepackage{graphicx}        
\usepackage{babel}



\begin{document}

\title{A Well-Balanced Scheme For Two-Fluid Flows In Variable Cross-Section
ducts}

\author{Philippe Helluy and Jonathan Jung}

\date{}

\maketitle
\begin{abstract}
We propose a finite volume scheme for computing two-fluid flows in
variable cross-section ducts. Our scheme satisfies a well-balanced
property. It is based on the VFRoe approach. The VFRoe variables are
the Riemann invariants of the stationnary wave and the cross-section.
In order to avoid spurious pressure oscillations, the well-balanced
approach is coupled with an ALE (Arbitrary Lagrangian Eulerian) technique
at the interface and a random sampling remap.
\end{abstract}

\section*{Introduction}

Classical finite volume solvers generally have a bad precision for
solving two-fluid interfaces or flows in varying cross-section
ducts. Several cures have been developed for improving the precision.
\begin{itemize}
\item For cross-section ducts, the well-balanced approach of Greenberg and
Leroux \cite{leroux96} (see also \cite{kroener2006} and \cite{jcam2011})
is an efficient tool to improve the precision.
\item For two-fluid flows the pressure oscillations phenomenon (see \cite{karni94}
and \cite{BHR04} for instance) can be cured by a recent tool developed
in \cite{CC2010} and \cite{BHMM2010}. It is based on an ALE (Arbitrary Lagrangian Eulerian)
scheme followed by a random sampling projection step.
\end{itemize}
In this paper, we show that is is possible to mix the two approaches
in order to design an efficient scheme for computing two-fluid flows
in variable cross-section ducts.

\section{A well-balanced two-fluid ALE solver}

\subsection{Model}
We consider the flow of a mixture of two compressible fluids (a gas
(1) and a liquid (2), for instance) in a cross-section duct. The time
variable is noted $t$ and the space variable along the duct is $x$.
We denote by $A(x)$ the cross-section at position $x$. The unknowns
are the density $\rho(x,t),$ the velocity $u(x,t)$, the internal
energy $e(x,t)$ and the fraction of gas $\varphi(x,t)$. Following
Greenberg and Leroux \cite{leroux96} it is now classical to consider
the cross-section $A$ as an artificial unknown. The equations are
the Euler equations in a duct, which read

\begin{eqnarray}
\partial_{t}(A\rho)+\partial_{x}(A\rho u) & = & 0,\label{first}\\
\partial_{t}(A\rho u)+\partial_{x}(A(\rho u^{2}+p)) & = & p\partial_{x}A,\\
\partial_{t}(A\rho E)+\partial_{x}(A(\rho E+p)u) & = & 0,\\
\partial_{t}(A\rho\varphi)+\partial_{x}(A\rho\varphi u) & = & 0,\\
\partial_{t}A & = & 0,\label{last}\end{eqnarray}
 with\begin{equation}
p=p(\rho,e,\varphi),\end{equation}
 \begin{equation}
E=e+\frac{u^{2}}{2}.\end{equation}
 Without loss of generality, in this paper we consider a stiffened
gas pressure law (see \cite{saurel1999} and included references) \begin{equation}
p(\rho,e,\varphi)=(\gamma(\varphi)-1)\rho e-\gamma(\varphi)\pi(\varphi).\end{equation}
 The mixture pressure law parameters $\gamma(\varphi)$ and $\pi(\varphi)$
are obtained from the pure fluid parameters $\gamma_{i}>1,\pi_{i}$,
$i=1,2$ thanks to the following interpolation, which is justified in \cite{BHR04}\begin{eqnarray}
\frac{1}{\gamma(\varphi)-1} & = & \varphi\frac{1}{\gamma_{1}-1}+(1-\varphi)\frac{1}{\gamma_{2}-1},\\
\frac{\gamma(\varphi)\pi(\varphi)}{\gamma(\varphi)-1} & =\varphi & \frac{\gamma_{1}\pi_{1}}{\gamma_{1}-1}+(1-\varphi)\frac{\gamma_{2}\pi_{2}}{\gamma_{2}-1}.\end{eqnarray}
 We define the vector of conservative variables\begin{equation}
W=(A\rho,A\rho u,A\rho E,A\rho\varphi,A)^{T}.\end{equation}
The conservative flux is\begin{equation}
F(W)=(A\rho u,A(\rho u^{2}+p),A(\rho E+p)u,A\rho\varphi u,0)^{T},\end{equation}
 and the non-conservative source term is\begin{equation}
S=(0,p\partial_{x}A,0,0,0),\end{equation}
 such that the system (\ref{first})-(\ref{last}) becomes
 \begin{equation}
\partial_{t}W+\partial_{x}F(W)=S(W).\label{systshort}\end{equation}
We define the vector of primitive variables\begin{equation}
Y=(\rho,u,p,\varphi,A)^{T}.\label{primitive}\end{equation}
 We define also the following quantities\begin{eqnarray}
Q & = & \text{mass flow rate}=\rho Au,\\
s & = & \text{entropy}=(p+\pi(\varphi))\rho^{-\gamma(\varphi)},\\
h & = & \text{enthalpy}=e+\frac{p}{\rho},\\
H & = & \text{total enthalpy}=h+\frac{u^{2}}{2}.\end{eqnarray}
 The entropy is solution of the partial differential equation\begin{equation}
Tds=de-\frac{p}{\rho^2}d\rho+\lambda d\varphi.\end{equation}
 It is useful to express also the pressure $p$ and the enthalpy $h$ as
functions of $(\rho,s,\varphi)$\begin{equation}
p=p(\rho,s,\varphi),\quad h=h(\rho,s,\varphi).\end{equation}
 Then in these variables the sound speed $c$ satisfies\begin{equation}
c^{2}=p_{\rho}=\rho h_{\rho}.\end{equation}
The jacobian matrix $F'(W)$ in system (\ref{systshort}) admits real eigenvalues\begin{equation}
\lambda_{0}=0,\quad\lambda_{1}=u-c,\quad\lambda_{2}=\lambda_{3}=u,\quad\lambda_{4}=u+c.\end{equation}
However, the system may be resonant (when $\lambda_0=\lambda_{1}$ or $\lambda_{0}=\lambda_{4}$.) The quantities $\varphi$, $s$, $Q$ and $H$ are independant Riemann invariants of the stationnary wave $\lambda_0$.
In the sequel, the vector of {}``stationary'' variables $Z$ will play a particular role \begin{equation}
Z=(A,\varphi,s,Q,H)^{T}.\label{stationary}\end{equation}
\subsection{VFRoe ALE numerical flux}
We recall now the principles of the VFRoe solver. We first consider
a arbitrary change of variables $U=U(W)$. In practice, we will take the set of
primitive variables $U=Y$ (\ref{primitive}) or the set of stationnary variables $U=Z$ (\ref{stationary}).
The vector $U$ satisfies a non-conservative set of equations\begin{equation}
\partial_{t}U+C(U)\partial_{x}U=0.\end{equation}
 The system (\ref{first})-(\ref{last}) is approximated by a finite volume scheme with cells $]x_{i-1/2},x_{i+1/2}[$, $i\in\mathbb{Z}$. We denote by $\tau$ the time step and by $\Delta x_i=x_{i+1/2}-x_{i-1/2}$ the size of cell $i$. We denote by
$W_{i}^{n}$ the conservative variables in cell $i$ at time step $n$.
The cross-section $A$ is approximated by a piecewise constant function,
$A=A_{i}$ in cell $i$. 

We consider first a very general scheme where the boundary of the cell $x_{i+1/2}$
moves at the velocity $v_{i+1/2}^{n}$ between time steps $n$ and $n+1$, thus we have\begin{equation}
x_{i+1/2}^{n+1}=x_{i+1/2}^{n}+\tau v_{i+1/2}^{n}.\end{equation}
In a VFRoe-type scheme, we have to define linearized
Riemann problems at interface $i+1/2$ between the state $W_{L}=W_{i}^{n}$
and $W_{R}=W_{i+1}^{n}$, we introduce\begin{equation}
\overline{U}=\frac{1}{2}(U_{L}+U_{R}).\end{equation}
 In this way, it is possibe to define\begin{equation}
\overline{W}=W(\overline{U}),\quad\overline{C}=C(\overline{U}).\end{equation}
 We then consider the linearized Riemann problem\begin{eqnarray}
\partial_{t}U+\overline{C}\partial_{x}U & = & 0,\\
U(x,\text{0)} & = & \left\{ \begin{array}{c}
U_{L}\text{ if }x<0,\\
U_{R}\text{ if }x>0.\end{array}\right.\end{eqnarray}
 We denote its solution by\begin{equation}
U(U_{L},U_{R},\frac{x}{t})=U(x,t).\end{equation}
Because of the stationary wave, $U(U_{L},U_{R},\frac{x}{t})$
is generally discontinuous at $x/t=0.$ We are then able to define
a discontinuous Arbitrary Lagrangian Eulerian (ALE) numerical flux
\begin{equation}
F(W_{L},W_{R},v^{\pm}):=F(W(U(U_{L},U_{R},v^{\pm})))-vW(U(U_{L},U_{R},v^{\pm})).\end{equation}
 The sizes of the cells evolve as\begin{equation}
\Delta x_{i}^{n+1}=\Delta x_{i}^{n}+\tau(v_{i+1/2}^{n}-v_{i-1/2}^{n}).\end{equation}
 If $v_{i+1/2}^{n}\leq0$ and $v_{i-1/2}^{n}\geq0$, the ALE
scheme is\begin{eqnarray}
\Delta x_{i}^{n+1}W_{i}^{n+1,-}-\Delta x_{i}^{n}W_{i}^{n}+\nonumber \\
\tau\left(F(W_{i}^{n},W_{i+1}^{n},v_{i+1/2}^{n,-})-F(W_{i-1}^{n},W_{i}^{n},v_{i-1/2}^{n,+})\right) & = & 0.\end{eqnarray}
 If $v_{i+1/2}^{n}>0$ then we have to add the following term to the
left of the previous equation\begin{equation}
\tau\left(F(W_{i}^{n},W_{i+1}^{n},0^{-})-F(W_{i}^{n},W_{i+1}^{n},0^{+})\right).\end{equation}
 If $v_{i-1/2}^{n}<0$ then we have to add also the following term\begin{equation}
\tau\left(F(W_{i-1}^{n},W_{i}^{n},0^{-})-F(W_{i-1}^{n},W_{i}^{n},0^{+})\right).\end{equation}

\subsection{ALE velocity}We have now to detail
the choice of the variable $U$ and the velocity $v$ according to
the data $W_{L}$ and $W_{R}$. The idea is to use the classical well-balanced
scheme everywhere but at the interface between the two fluids, where we use the Lagrange flux. When our initial data satisfy 
$\varphi \in \{0,1\}$, the algorithm reads
\begin{itemize}
\item If we are not at the interface, i.e. if $\varphi_{L}=\varphi_{R}$, we  take $U=Z$ and $v=0$. This choice corresponds
to the VFRoe well-balanced scheme described in \cite{jcam2011}.

\item If we are at the interface, i.e. if $\varphi_{L}\neq\varphi_{R}$
then we choose $U=Y$. This choice ensures that the linearized Riemann
solver presents no jump of pressure and velocity
at the contact discontinuity. We thus denote by $u^{*}(W_{L},W_{R})$
and $p^{*}(W_{L},W_{R})$ the velocity and the pressure at the contact.
We take $v=u^{*}(W_{L},W_{R})$, $A^{*}=A_{L}$ if $v<0$ and $A^{*}=A_{R}$
if $v>0$. The lagrangian numerical flux then takes the form\begin{equation}
F(W_{L},W_{R},v^{\pm})=(0,A^{*}p^{*},A^{*}u^{*}p^{*},0,-A^{*}u^*)^T.\end{equation}
\end{itemize}

\subsection{Glimm remap}

We go back to the original Euler grid by the Glimm procedure.

We construct a sequence of pseudo-random numbers $\omega_{n}\in[0,1[.$
In practice, we consider the $(5,3)$ van der Corput sequence \cite{BHMM2010}.
According to this number we take\begin{equation}
W_{i}^{n+1}=W_{i-1}^{n+1,-}\text{ if }\omega_{n}<\frac{\tau_{n}}{\Delta x_{i}}\max(v_{i-1/2}^{n},0),
\end{equation}
\begin{equation}
W_{i}^{n+1}=W_{i+1}^{n+1,-}\text{ if }\omega_{n}>1+\frac{\tau_{n}}{\Delta x_{i}}\min(v_{i+1/2}^{n},0),
\end{equation}
\begin{equation}
W_{i}^{n}=W_{i}^{n+1,-}\text{ if }\frac{\tau_{n}}{\Delta x_{i}}\max(v_{i-1/2}^{n},0)\leq\omega_{n}\leq1+\frac{\tau_{n}}{\Delta x_{i}}\min(v_{i+1/2}^{n},0).\end{equation}

\subsection{Properties of the scheme}
The constructed scheme has many interesting properties:
\begin{itemize}
\item it is well-balanced in the sense that it preserves exactly all stationary states (i.e. initial data for which the quantities $\varphi, s,Q,H$ are constant);
\item for constant cross-section ducts, it computes exactly the contact discontinuities, with no smearing of the density and the mass fraction;
\item if at the initial time the mass fraction is in $\{0,1\}$, then this property is exactly preserved at any time.
\end{itemize}
For detailed proofs, we refer to \cite{jcam2011} and \cite{BHMM2010}. Some other subtleties are given in the same references. For instance, the change of variables $Z=Z(W)$ is not always invertible. This implies to define a special procedure for constructing completely rigorously the well-balanced VFRoe solver.

\section{Numerical results}

In order to test our algorithm, we consider a Riemann problem for which we know the exact solution. The initial data are discontinuous at $x=1$. The data of the problem are given in Table \ref{table}

\begin{table}
\begin{center}
\begin{tabular}{|c|c|c|}
\hline 
quantity & Left & Right\tabularnewline
\hline
\hline 
$\rho$ & 2 & 3.230672602\tabularnewline
\hline 
$u$ & 0.5 & -0.4442565900\tabularnewline
\hline 
$p$ & 1 & 12\tabularnewline
\hline 
$\varphi$ & 1 & 0\tabularnewline
\hline 
$A$ & 1.5 & 1\tabularnewline
\hline
\end{tabular}
\caption{Numerical results. Data of the Riemann problem\label{table}}
\end{center}
\end{table}

\begin{figure}
\begin{center}
\includegraphics[angle=-90,width=8cm]{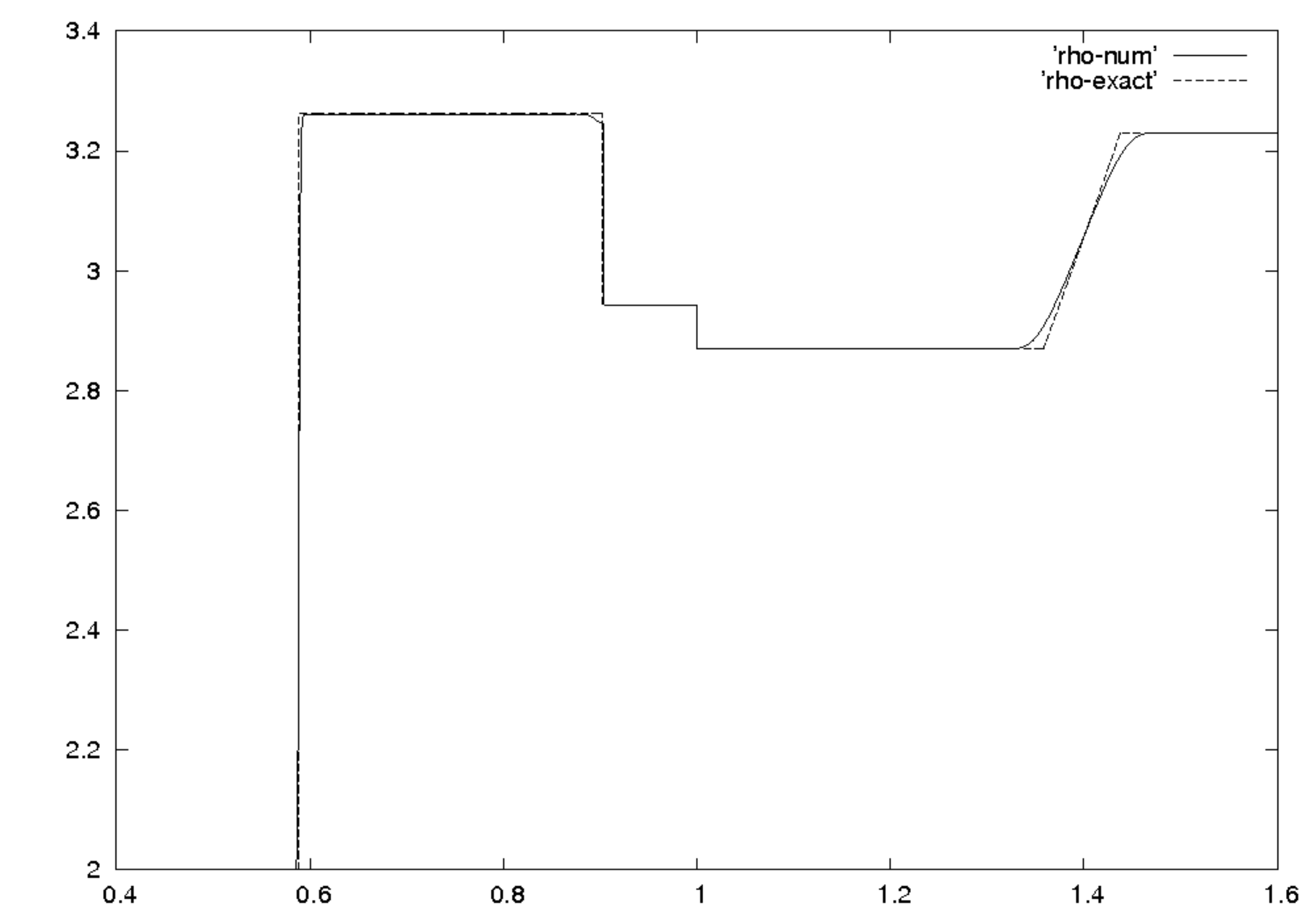}
\end{center}
\caption{Two-fluid, discontinuous cross-section Riemann problem. Density plot.
Comparison of the exact solution (dotted line) and the approximate
one (continuous line)\label{fig:rho}}

\end{figure}

\begin{figure}
\begin{center}
\includegraphics[angle=-90,width=8cm]{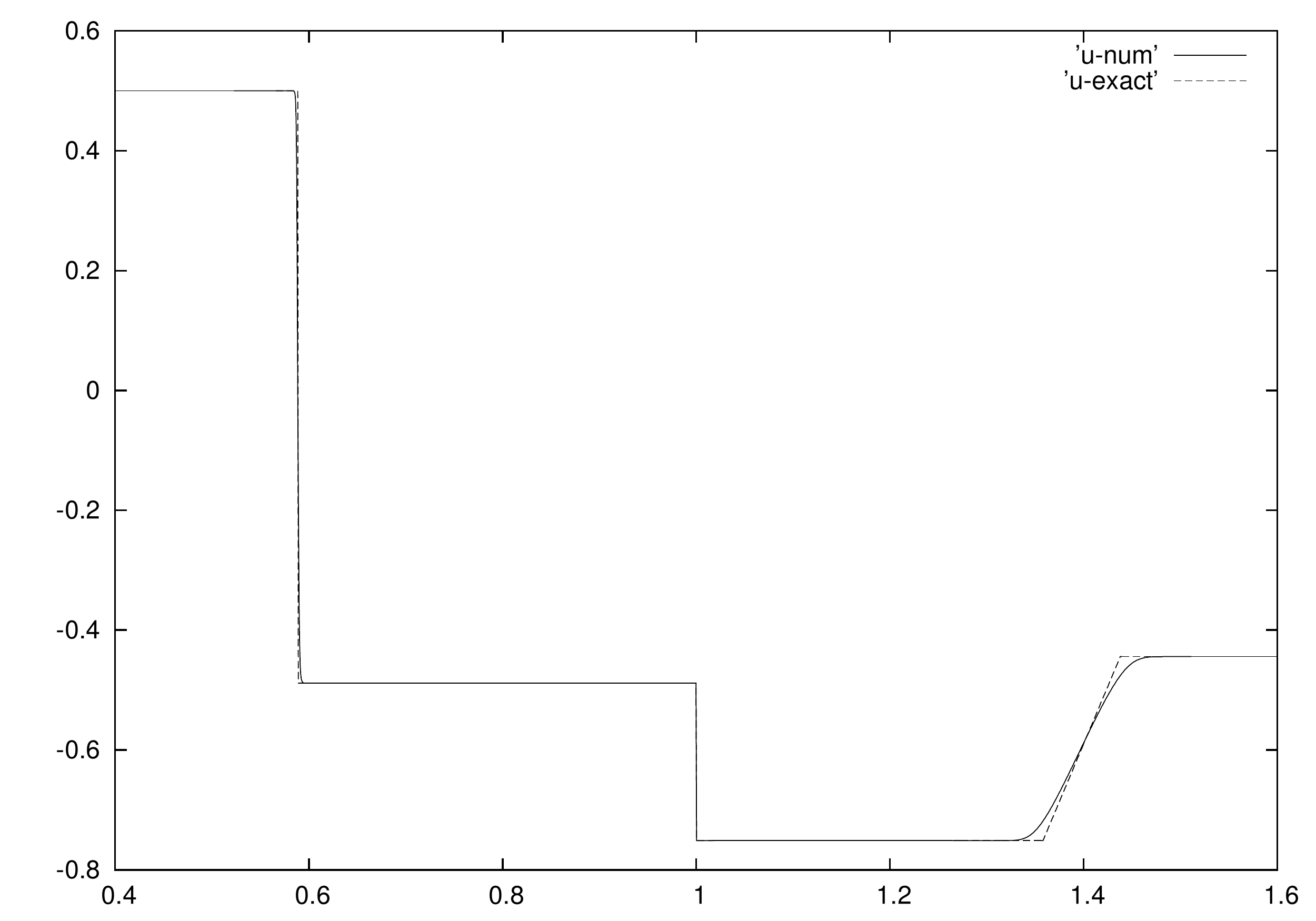}
\end{center}
\caption{Two-fluid, discontinuous cross-section Riemann problem. Pressure plot.
Comparison of the exact solution (dotted line) and the approximate
one (continuous line)\label{fig:u}}

\end{figure}

\begin{figure}
\begin{center}
\includegraphics[angle=-90,width=8cm]{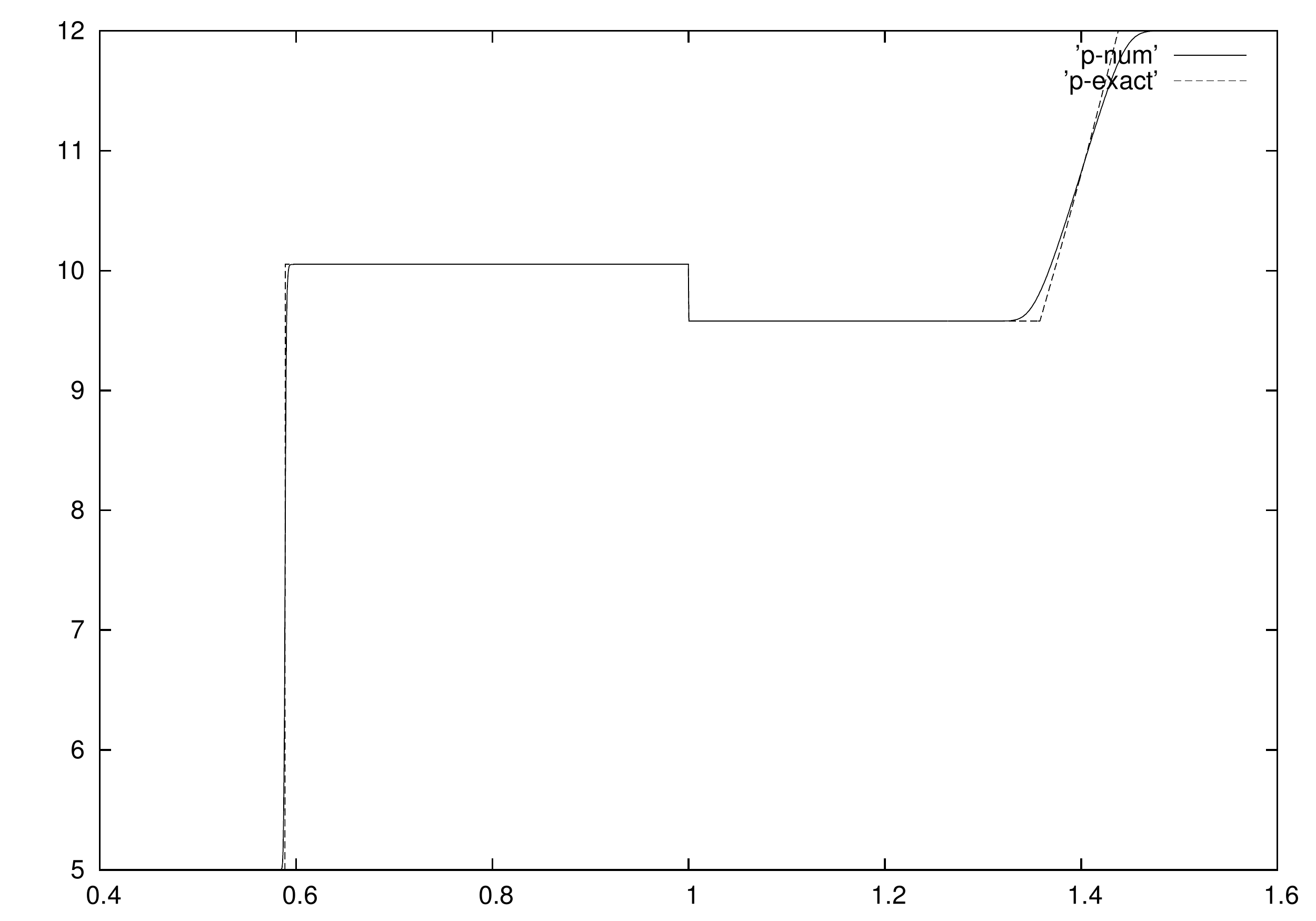}
\end{center}
\caption{Two-fluid, discontinuous cross-section Riemann problem. Velocity plot.
Comparison of the exact solution (dotted line) and the approximate
one (continuous line)\label{fig:p}}

\end{figure}

The pressure law parameters are $\gamma_{1}=1.4$, $\pi_{1}=0$, $\gamma_{2}=1.6$
and $\pi_{2}=2.$ We compute the solution on the domain $[0.4;1.6]$
with approximately $2000$ cells. The final time is $T=0.2$ and the
CFL number is $0.6$. The density, the velocity and the pressure
are represented on Figures \ref{fig:rho}, \ref{fig:u} and \ref{fig:p}.
We observe an excellent agreement between the exact and the approximate
solution. The mass fraction is not represented: it is not smeared at all
and perfectly matches the exact solution.


\section{Conclusion}

We have constructed and validated a new scheme for computing two-fluid
flows in variable cross-section ducts. Our scheme relies on two ingredients:
\begin{itemize}
\item a well-balanced approach for dealing with the varying cross-section;
\item a Lagrange plus remap technique in order to avoid pressure oscillations
at the interface. The random sampling remap ensures that the interface
is not diffused at all.
\end{itemize}
On preliminary test cases, our approach gives very satisfactory results.
We intend to apply it to the computation of the oscillations of cavitation
bubbles. More results will be presented at the conference.

The authors wish to thank Jean-Marc H\'erard for many fruitful discussions.

\end{document}